\documentclass[reqno,12pt]{amsproc}
\newtheorem{theo}{Теорема}

\newtheorem{lem}{Лемма}

\newtheorem{df}{Определение}

\newcommand\eps\varepsilon
\newcommand \bs\boldsymbol
\newcommand\ph\varphi
\newcommand\kap\varkappa

\newcommand\conv{\mathrm{conv}\,}

\usepackage{enumerate}
\usepackage{graphicx}
\usepackage{float}
\usepackage[cp1251]{inputenc}
\usepackage[english, russian]{babel}

\begin{document}\title[О периодических движениях]
{О периодических движениях в одной задаче с сухим трением}

\author[]{Олег Зубелевич \\ \\\tt
 Кафедра теоретической механики и мехатроники,  \\
Механико-математический факультет,\\
МГУ им. М.В. Ломоносова, Москва\\
\\
Математический институт им. В.А. Стеклова Российской академии наук
 }

\date{}
\thanks{Исследование выполнено за счет гранта Российского научного фонда
(проект № 19-71-30012).}
\subjclass[2000]{  65L05, 34A60, 34C25, 37N05, 70F40, 74H45}
\keywords{Differential inclusions, Coulomb friction, periodic solutions, nonsmooth dynamics.}

\begin{abstract}
\textsc{ Abstract.} We consider a point mass  on a horizontal plane. The motion of the plane is given. The plane moves periodically such that all its points  have  congruent closed trajectories. 
There is the Coulomb friction between the point mass and the plane. 
We show that there exists a motion  of the point mass such that the velocity of the point mass is an absolutely continuous periodic function. 
Mathematically the result is expressed as an existence  theorem for a periodic solution to some differential inclusion. 

-----

В статье рассмотрено движение частицы по горизонтальной вибрирующей  плоскости с кулоновским  трением. Доказана теорема существования решений с периодической по времени скоростью.

\end{abstract}

\maketitle
\section{Введение, постановка задачи}
Вопросы феноменологии и математического моделирования  систем с сухим трением до сих пор не до конца  прояснены, см., например, \cite{Kozlov}, \cite{Zhuravlev}. Известно, что в системах  с сухим трением задача Коши может быть некорректной, а сама гипотеза сухого трения может приводить к парадоксам -- результатам, находящимся в явном противоречии с физикой
 \cite{Penleve}. 
 
 По мнению автора, важным вкладом в теорию сухого трения является работа А. Ф. Филиппова \cite{fil}, результаты которой позволяют, по крайней мере в ряде случаев, построить корректную модель системы с сухим трением и от феноменологии перейти к динамике, т. е.  к изучению поведения решений соответствующей системы дифференциальных уравнений. 
 
 В этой заметке мы рассмотрим приложения результатов А. Ф. Филиппова к следующей модельной задаче.

{\it Шероховатая горизонтальная плоскость $\Pi$ перемещается  поступательно по заданному закону,  не меняя высоты, таким образом, что каждая ее точка совершает периодическое движение с периодом $\omega>0$.

По плоскости в стандартном поле силы тяжести  скользит частица $P$ массы $m$, на частицу действует сила кулоновского  трения с коэффициентом $\sigma>0$,  а также сила тяжести, которая уравновешена нормальной реакцией плоскости.

Спрашивается, существуют ли в данной системе периодические движения?}

По-видимому, ответ на этот вопрос непрост.  Однако, мы докажем, что существуют решения с периодической по времени скоростью.

  Удобно выбрать размерности физических величин так, что $\omega=1,\quad g=1,\quad m=1$. 

По условию, скорости всех точек плоскости $\Pi$ относительно лабораторной системы отсчета одинаковы. Обозначим эту скорость за $\bs V(t)$. Вектор-функция $\bs V\in C^1([0,\infty),\mathbb{R}^2)$ считается заданной; ее период равен $1$.  Далее будем считать, что $\|\bs{\dot V}\|_{C[0,1]}\ne 0$ -- в противном случае задача тривиальна.   Через $\bs u$ обозначим скорость частицы относительно плоскости $\Pi$.

Запишем уравнения движения:
\begin{equation}\label{dfff6}\boldsymbol {\dot u}=\bs F(t,\bs u),\quad \bs F(t,\bs u)=-\sigma\frac{\bs u}{|\bs u|}-\bs {\dot V},\quad \bs u\in\mathbb{R}^2.\end{equation} Через $|\cdot|$ мы обозначаем стандартную евклидову норму.

Правая часть  уравнения (\ref{dfff6}) не может быть продолжена на прямую $\{t\in\mathbb{R},\,\bs u=0\}$ расширенного фазового пространства  непрерывно. Таким образом, обычное определение решения не имеет смысла. Мы будем искать периодические решения в классе обобщенных решений в смысле Филиппова \cite {fil}.

\subsection{Обобщенные решения в смысле Филиппова} \label{sdfggggg}Векторы нашего физического евклидового пространства мы, как и выше, обозначаем жирным шрифтом. Векторы абстрактного $m-$мерного пространства мы обозначаем обычным шрифтом:
$x=(x^1,\ldots,x^m)^T\in\mathbb{R}^m$.

Через $\conv{U}$ мы обозначаем замкнутую выпуклую оболочку множества $U\subset\mathbb{R}^m$. Напомним, что замкнутой выпуклой оболочкой называется пересечение всех выпуклых замкнутых множеств, содержащих $U$.

Через
$$B_s(x)=\{y\in\mathbb{R}^m\mid |y-x|<s\}$$ обозначим открытый шар в $\mathbb{R}^m$;
через $\mu$ обозначим меру Лебега в $\mathbb{R}^m$.

Пусть  $D\subset\mathbb{R}^m$ -- некоторая область.
Отображение $$f:[t_1,t_2]\times D\to \mathbb{R}^m$$  измеримо по Лебегу.

\begin{df}[\cite{fil}]\label{xdfggg}Абсолютно непрерывную функцию $x:[t_1,t_2]\to D$ назовем обобщенным решением  уравнения 
\begin{equation}\label{ffnggg}\dot x=f(t,x),\end{equation}
если для почти всех $t\in[t_1,t_2]$ имеет место включение
\begin{equation}\label{sfgggg}
\dot x(t)\in\bigcap_{\delta>0}\bigcap_{N}\conv f(t,B_\delta(x(t))\backslash N).   \end{equation}
Пересечение берется по всем измеримым множествам $N\subset\mathbb{R}^m$ нулевой меры Лебега: $\mu(N)=0$.

Через $f(t,B_\delta(x(t))\backslash N)\subset\mathbb{R}^m$ обозначен образ множества $B_\delta(x(t))\backslash N$ при отображении $x\mapsto f(t,x)$ с фиксированным $t$.
\end{df}
Если векторное поле $f$ непрерывно, то множество в правой части формулы (\ref{sfgggg}) состоит из единственного элемента $f(t,x(t))$, и определение \ref{xdfggg} переходит в стандартное определение решения дифференциального уравнения.

В соответствии с теоремой 2 \cite{fil}, свойство функции быть обобщенным решением не зависит от выбора системы координат  $x$.

\section{Основная теорема}
Как уже отмечалось, функция $\bs F$ не может быть продолжена на ось $\{(t,0)\}$ непрерывно. 
Однако то, как именно мы продолжим функцию $\bs F$ на прямую $\{(t,0)\}$, не имеет значения с точки зрения определения \ref {xdfggg}. При этом данное определение дает именно тот результат, который мы и ожидаем получить из теории кулоновского трения.

Действительно, множество 
\begin{equation}W=\label{dfh66}\bigcap_{\delta>0}\,\bigcap_{\mu(N)=0}\conv \Big(\bs F(t,B_\delta(0)\backslash N)\Big)=\overline B_\sigma(-\bs {\dot V}(t))\end{equation}
(чертой обозначено замыкание шара) не зависит от значений $\bs F(t,0)$, т. к. одно из множеств нулевой меры $N$ обязательно содержит точку $\bs u=0$, и значение $\bs F(t,0)$, каким бы мы его ни выбрали, на множество (\ref{dfh66}) не влияет.

Если на интервале времени $(t',t'')$ скорость $\bs u(t)$ не равна нулю, то $\bs u(t)$ -- это гладкая функция времени, удовлетворяющая уравнению (\ref{dfff6}) в обычном смысле. 

Если частица  покоится  относительно плоскости $\Pi$ на  интервале  $t\in (\tau',\tau'')$,  то, в соответствии с   определением обобщенного решения, имеем: 
\begin{equation}\label{srtg5t}\bs {\dot u}(t)=0\in W,\quad t\in (\tau',\tau'').\end{equation}
 Тогда из формулы  (\ref{srtg5t}) следует,  что   $$|\bs {\dot V}(t)|\le\sigma.$$
 Это согласуется с определением кулоновского трения.

Заметим, что по теореме 10 \cite{fil}, обобщенное решение системы  (\ref{dfff6}) единственно вперед при любом заданном начальном условии $\bs u(t_0)$.
Ключевую роль при проверке условий теоремы единственности играет то, что сила трения
$$\bs\Phi(\bs u)=-\sigma\frac{\bs u}{|\bs u|}$$
монотонна:
$$(\bs v-\bs w,\bs\Phi(\bs v)-\bs\Phi(\bs w))\le 0,\quad\forall \bs v,\bs w\ne 0.$$

\begin{theo} \label{xdfggg}Система (\ref{dfff6}) имеет $1$- периодическое обобщенное решение $\bs u(t).$ \end{theo}

Полученное в теореме \ref{xdfggg} решение соответствует движению, при котором относительная скорость частицы $P$ может обращаться в ноль, частица может некоторое время не менять своего положения относительно плоскости $\Pi$, а потом опять начать движение.

\section{Доказательство теоремы \ref{xdfggg}}
\subsection{Периодические решения аппроксимирующей задачи}
Рассмотрим последовательность систем
\begin{equation}\label{dfggg}
\boldsymbol {\dot u}=\bs F(t,\bs u)-\frac{1}{k}\bs u,\quad k\in\mathbb{N}.\end{equation}
\begin{lem}\label{dsfgffda4}
Предположим, что все системы  (\ref{dfggg}) имеют $1-$периодические обобщенные решения $\bs u_k(t)$, определенные на $[0,\infty)$. Тогда для всех $k$ верна оценка
\begin{equation}\label{xdfhhhh}\|\bs u_k\|_{C[0,1]}\le \rho.\end{equation}
Положительная постоянная $\rho$ не зависит от $k$.\end{lem}
Отметим, что в области $\{\bs u\ne 0\}$ система (\ref{dfggg}) является гладкой, а ее обобщенные решения совпадают с классическими, и для них верны все теоремы теории дифференциальных уравнений с гладкой правой частью.

\subsubsection{Доказательство леммы \ref{dsfgffda4}}

Доказательство разобьем на два случая. 

Первый случай: периодическое решение $\bs u_k$ обращается в ноль при некоторых $t$. 

Зафиксируем $k$.
Пусть $t^*$ -- точка максимума функции $|\bs u_k(t)|$. 
Будем считать, что $\bs u_k(t^*)\ne 0$,  иначе доказывать нечего.
Тогда, очевидно, найдутся точки $$t^-<t^+,\quad t^+-t^-\le 1,\quad \bs u_k(t^\pm)=0$$ и последовательности $$t_j^{\pm}\to t^{\pm},\quad t^-<t_{j+1}^-<t_j^-<t^*<t_j^+<t_{j+1}^+<t^+$$ такие, что
$$\min_{t\in[t_j^-,t_j^+]}|\bs u_k(t)|>0,\quad|\bs u_k(t_j^+)|=|\bs u_k(t_j^-)|\to 0,\quad j\to \infty. $$

На интервале $[t_j^-,t_j^+]$ функция $\bs u_k$ непрерывно дифференцируема и является решением системы (\ref{dfggg}) в стандартном смысле. 

Домножим скалярно левую и правую часть уравнения (\ref{dfggg}) на $\bs{\dot u}_k$:
\begin{equation}\label{dfghhh}|\bs{\dot u}_k|^2=-\sigma\frac{d}{dt}|\bs{ u}_k|-(\bs{\dot V},\bs{\dot u}_k)-\frac{1}{2k}\frac{d}{dt}|\bs{ u}_k|^2.\end{equation}
Интегрируя это равенство по отрезку $[t_j^-,t_j^+]$,  находим:
$$\|\bs{\dot u}_k\|_{L^2(t_j^-,t_j^+)}^2=-\int_{t_j^-}^{t_j^+}(\bs {\dot V},\bs{\dot u}_k)dt\le\|\bs {\dot V}\|_{L^2(t_j^-,t_j^+)}
\|\bs{\dot u}_k\|_{L^2(t_j^-,t_j^+)}.$$
Переходя к пределу при $j\to \infty$, по теореме Леви получаем
$$\|\bs{\dot u}_k\|_{L^2(t^-,t^+)}\le \|\bs {\dot V}\|_{L^2(t^-,t^+)}.$$
Из формулы
$$\bs u_k(t^*)=\int_{t^-}^{t^*}\bs {\dot u}_k(t)dt$$ вытекает оценка
$$|\bs u_k(t^*)|\le \sqrt{t^*-t^-}\|\bs{\dot u}_k\|_{L^2(t^-,t^*)}\le \|\bs {\dot V}\|_{L^2(0,1)}. $$
Что и доказывает лемму для первого случая.

Переходим ко второму случаю: $|\bs u_k(t)|>0$ при всех  $t$. В этом случае $\bs u_k(t)$ -- непрерывно дифференцируемая на $[0,1]$  функция, которая является классическим решением уравнения (\ref{dfggg}).

Интегрируя формулу (\ref{dfghhh}) по отрезку $[0,1]$, мы получаем оценку
\begin{equation}\label{frh90}
\|\bs{\dot u}_k\|_{L^2(0,1)}\le \|\bs {\dot V}\|_{L^2(0,1)}.\end{equation}
Домножим скалярно левую и правую часть уравнения (\ref{dfggg}) на $\bs{ u}_k$:
$$\frac{1}{2}\frac{d}{dt}|\bs u_k|^2=-\sigma|\bs u_k|-(\bs{\dot V},\bs u_k)-\frac{1}{k}|\bs u_k|^2.$$
Интегрируя это равенство по отрезку $[0,1]$ и применяя формулу интегрирования по частям, имеем:
$$\sigma\|\bs u_k\|_{L^1(0,1)}=\int_0^1(\bs{ V},\bs{\dot u}_k)dt-\frac{1}{k}\|\bs u_k\|_{L^2(0,1)}^2.$$
Откуда
\begin{equation}\label{gfh667}
\sigma\|\bs u_k\|_{L^1(0,1)}\le  \|\bs { V}\|_{L^2(0,1)} \|\bs {\dot u}_k\|_{L^2(0,1)}.\end{equation}
Из формул (\ref{frh90}), (\ref{gfh667}) вытекает, что последовательность $\{\bs u_k\}$ ограничена в пространстве Соболева:
$$\sup_k\|\bs u_k\|_{H^{1,1}(0,1)}<\infty.$$
Теперь утверждение леммы вытекает из теоремы вложения: $H^{1,1}(0,1)\subset L^\infty(0,1)$ \cite{adams}.

Лемма доказана.

Из результатов работы \cite{fil} вытекает, что для каждого начального условия $\bs u_0=\bs u(t_0)$ система (\ref{dfggg}) имеет на некотором интервале  $t\in[ t_0,t_1)$ обобщенное решение, это решение единственно вперед и непрерывно зависит от $\bs u_0$.

Более того, рассуждая как и в случае гладких систем, можно показать, что если обобщенное решение $\bs u(t)$ системы (\ref{dfggg})  определено и ограничено на интервале $[t_0,t_1),$ то оно продолжается вправо, за точку $t_1$.

Пусть  $\bs u(t)$ -- какое-нибудь решение системы (\ref{dfggg}). Покажем, что если  
$$\bs u(t_0)\in \overline B_{R_k}(0),\quad R_k=k\|\bs{\dot V}\|_{C[0,1]},$$ 
то $\bs u(t)\in  B_{R_k}(0),\quad \forall t>t_0$.

Действительно, пока решение $\bs u$ находится рядом с границей шара, оно является классическим.
Домножим левую и правую часть  (\ref{dfggg}) на это решение:
$$
\frac{1}{2}\frac{d}{dt}|\bs u|^2=-\sigma|\bs u|-(\bs{\dot V},\bs{u})-\frac{1}{k}|\bs u|^2.$$
Откуда
$$\frac{1}{2}\frac{d}{dt}|\bs u|^2\le (-\sigma+|\bs{\dot V}|)|\bs u|-\frac{1}{k}|\bs u|^2.$$
Из последней формулы следует, что если решение $\bs u$ находится рядом с границей шара, то 
$$\frac{d}{dt}|\bs u|^2<0,$$ и решение отходит от границы.

Таким образом, для каждой системы (\ref{dfggg}) определено непрерывное отображение Пуанкаре $T_k:\overline B_{R_k}(0)\to \overline B_{R_k}(0)$. По теореме Брауэра \cite{Loran},  это отображение имеет неподвижную точку, которая является начальным условием для $1-$периодического обобщенного решения системы (\ref{dfggg}).

Как и выше, будем обозначать эти обобщенные периодические решения $\bs u_k(t)$. Они удовлетворяют оценке (\ref{xdfhhhh}).

\subsection{Теорема об аппроксимации}
 Пусть в добавление к предположениям раздела \ref{sdfggggg}, отображение $f$ ограничено почти всюду:
$$\|f\|_{L^\infty([t_1,t_2]\times D)}<\infty.$$

\begin{theo}[\cite{fil}]\label{df568}
 Предположим, что существует последовательность  измеримых функций 
 $$g_k,f_k:[t_1,t_2]\times D\to \mathbb{R}^m,\quad k\in\mathbb{N}$$
 таких, что для почти всех $(t,x)$ имеет место включение
$$f_k(t,x)\in \bigcap_{\mu(N)=0}\conv f(t,B_{\delta_k}(x)\backslash N),\quad \delta_k\to 0,\quad k\in\mathbb{N};$$ 
 и найдется последовательность функций $q_k(t)\in L^1(t_1,t_2)$ таких, что для почти всех $(t,x)$
 верны оценки
 $$|g_k(t,x)|\le q_k(t),\quad \int_{t_1}^{t_2}q_k(t)dt\to 0,\quad k\to\infty.$$

Тогда если последовательность обобщенных решений $x_k:[t_1,t_2]\to D$
систем
\begin{equation}\label{xdfb43}
\dot x_k=f_k(t,x_k)+g_k(t,x_k)\end{equation} равномерно ограничена:
$$\sup_k\|x_k\|_{C[t_1,t_2]}<\infty,$$
то она содержит равномерно сходящуюся на $[t_1,t_2]$ подпоследовательность $\{x_{k_j}\}$, и предел этой подпоследовательности  является обобщенным решением системы (\ref{ffnggg}).
\end{theo}

\subsection{Периодические решения системы (\ref{dfff6})}
Теперь утверждение теоремы \ref{xdfggg} вытекает непосредственно из теоремы \ref{df568}. Действительно,
в качестве множества $D$ возьмем шар $B_{\rho+1}(0)$.  В качестве системы (\ref{xdfb43}) возьмем систему (\ref{dfggg}),
причем $g_k$ это $-\bs  u/k$.
Ясно, что
$$\Big|\frac{1}{k}\bs  u\Big|\le\frac{\rho+1}{k}=q_k.$$
В качестве решений $x_k$ берем решения $\bs u_k$. В качестве системы (\ref{ffnggg}) -- систему (\ref{dfff6}).

Теорема доказана.

\end{document}